\documentclass[11pt]{article}
\usepackage{cite}
\usepackage{mathrsfs}
\usepackage{amssymb,amsfonts,amsmath}
\usepackage{enumerate}
\usepackage{indentfirst}
\usepackage{latexsym,bm}
\usepackage[noend]{algpseudocode}
\usepackage{algorithmicx,algorithm}
\usepackage{algorithm}
\usepackage{makeidx}
\usepackage{fancybox}
\usepackage{color}
\usepackage{multicol}
\usepackage[latin1]{inputenc}
\usepackage{graphicx}
\usepackage{epstopdf}
\usepackage{pstricks,pst-node,pst-text,pst-3d}
\usepackage{tikz}
\newtheorem{thm}{Theorem}[section]
\newtheorem{cor}[thm]{Corollary}
\newtheorem{lem}[thm]{Lemma}
\newtheorem{op}{Problem}[section]

\newtheorem{conj}{Conjecture}[section]

\newenvironment{pf}{{\noindent \it \bf Proof:}}{{\hfill$\Box$}\\}

\begin{document}

\title{\bf Compositions of Digraphs: A Survey}
\author{Yuefang Sun\\
School of Mathematics and Statistics, Ningbo University,\\
Zhejiang 315211, P. R. China, \\sunyuefang@nbu.edu.cn}
\maketitle

\tableofcontents
\newpage

\begin{abstract}
In this survey we overview known results and get several new results on  digraph compositions which generalize several classes of digraphs, such as quasi-transitive digraphs. After an introductory section, the paper is divided into six sections: connectivity and linkages, kings and kernels, paths, cycles, acyclic spanning subdigraphs, strong spanning subdigraphs. This survey also contains some conjectures and open problems for further study.
\vspace{0.3cm}\\
{\bf Keywords:} digraph composition; semicomplete composition; transitive composition; connectivity; linkage; king; kernel; Hamiltonicity; path-cycle subdigraph; Path Partition Conjecture; pancyclicity; acyclic spanning subgraph; branching; strong spanning subdigraph.
\vspace{0.3cm}\\ {\bf AMS subject
classification (2020)}: 05C12, 05C20, 05C38, 05C40, 05C45, 05C69, 05C70, 05C75, 05C76, 05C85.
\end{abstract}

\section{Introduction}\label{sec:intro}

In this section, we give the motivation for introducing the notion of digraph compositions by discussing the relationship between digraph compositions and several other digraph classes. Some notation and terminology are also given.

\subsection{Digraph compositions and related digraph classes}

We refer the readers to \cite{Bang-Jensen-Gutin, Bondy} for graph-theoretical notation and terminology not given here. Throughout this survey, unless otherwise stated, paths and cycles are always assumed to be directed, and all digraphs considered in this paper have no parallel arcs
or loops. We use $[n]$ to denote the set of all natural numbers from 1 to $n$.

A digraph $D$ is {\em connected} if $U(D)$ is connected, where $U(D)$ is the {\em underlying graph} of $D$. A digraph is {\em acyclic} if it has no directed cycle.
A digraph $D$ is {\em semicomplete} if for every pair $x,y$ of distinct vertices of $D$, there is at least one arc between $x$ and $y.$ In particular, a {\em tournament} is a semicomplete digraph without a 2-cycle. 
A digraph $D$ is {\em transitive} (resp. {\em quasi-transitive}), if for any triple $x,y,z$ of distinct vertices of $D$, if $xy$ and $yz$ are arcs of $D$ then $xz\in A(D)$ (resp. then either $xz$ or $zx$ or both are arcs  of $D$).

We often use the following operation called ``composition'' to construct ``bigger'' digraphs from ``smaller'' ones.
Let $T$ be a digraph with vertices $u_1, \dots, u_t$ ($t\ge 2$) and let $H_1, \dots, H_t$ be digraphs such that $H_i$ has vertices $u_{i,j_i}, j_i\in [n_i].$  Let $n_0=\min\{n_i\mid i\in [t]\}$. Then the {\em composition} $Q=T[H_1, \dots, H_t]$ is a digraph with vertex set $$V(Q)=\{u_{i,j_i}\mid i\in [t], j_i\in [n_i]\}$$ and arc set \[
\left(\bigcup^t_{i=1}A(H_i) \right) \bigcup  \left( \bigcup_{u_iu_p\in A(T)} \{u_{ij_i}u_{pq_p} \mid j_i\in [n_i], q_p\in [n_p]\} \right).
\] The composition $Q=T[H_1, \dots, H_t]$ is {\em proper} if $Q\neq T$, i.e. at least one $H_i$ is nontrivial, and is {\em semicomplete} (resp. {\em transitive}) if $T$ is semicomplete (resp. transitive).

Digraph composition is related to another operation. Let $T$ be a digraph with vertex set $V(T)=\{u_i\mid i\in [t]\}$ and let $K$ be another digraph. By {\em blowing up $u_i$ into $K$ in $T$} we mean the operation that substitutes the digraph $K$ for the vertex $u_i$ in $T$, that is, produces a new digraph $D=[\{u_1\}, \dots, \{u_{i-1}\}, K, \{u_{i+1}\}, \dots, \{u_t\}]$. Clearly, the digraph composition $Q=T[H_1, \dots, H_t]$ can be obtained from $T$ and $\{H_i\mid i\in [t]\}$ by blowing up $u_i$ into $H_i$ for each $i\in [t]$.

The notion of digraph composition can be used to define a class of decomposable digraphs. Note that if $Q=T[H_1, \dots, H_t]$, then $T, H_1, \dots, H_t$ are induced subdigraphs of $Q$ and we say that $Q$ is {\em decomposable} (into $T, H_1, \dots, H_t$). Let $\Phi$ be a class of digraphs. A digraph $Q$ is {\em $\Phi$-decomposable} if $Q\in \Phi$ or $Q=T[H_1, \dots, H_t]$ for some $T\in \Phi$ with $t=|T|\geq 2$ and some choice of digraphs $H_1, \dots, H_t$, and this decomposition is called a {\em $\Phi$-decomposition}.
In particular, we say that a digraph $Q$ is {\em totally $\Phi$-decomposable} if either $Q\in \Phi$ or $Q=T[H_1, \dots, H_t]$ with $T\in \Phi$ and $H_i$ totally $\Phi$-decomposable for $i\in [t]$. The {\em total $\Phi$-decomposition} of $Q$ is inductively defined as the sequence:
\[
\left\{
   \begin{array}{ll}
     Q, &\mbox {if $Q\in \Phi$;}\\
     T, L_1, \dots, L_t,~$where$~L_i~$is~the~total$~\Phi-$decomposition~of~$H_i, &\mbox {Otherwise.}
   \end{array}
   \right.
\]

The {\em first flyer} of the total $\Phi$-decomposition of $Q=T[H_1, \dots, H_t]$, namely, $T, H_1, \dots, H_t$ is called the {\em $\Phi$-decomposition} of $D$. 

There are several important digraph classes which are subclasses of digraph compositions. Interestingly, some of them are also subclasses of totally $\Phi$-decomposable digraphs.
If $Q=T[H_1, \dots, H_t]$ and none of the digraphs $H_1, \dots, H_t$ has an arc, then $Q$ is an {\em extension} of $T$. For any set $\Phi$ of digraphs, $\Phi^{ext}$ denotes the (infinite) set of all extensions of digraphs in $\Phi$, which are called {\em extended $\Phi$-digraphs}. For example, if $\Phi$ is the set of all semicomplete digraphs, then $\Phi^{ext}$ denotes the set of all extended semicomplete digraphs. Clearly, extended semicomplete digraphs are totally $\Phi_1$-decomposable, where $\Phi_1$ is the set of all semicomplete digraphs and acyclic digraphs. 

The following theorem by Bang-Jensen and Huang gives a complete characterization of quasi-transitive digraphs which shows that quasi-transitive digraphs are totally $\Phi_1$-decomposable. The decomposition below is called the {\em canonical decomposition} of a quasi-transitive digraph.

\begin{thm}\label{intro01}\cite{Bang-Jensen-Huang1}
Let $D$ be a quasi-transitive digraph. Then the following assertions hold:
\begin{description}
\item[(a)] If $D$ is not strong, then there exists a transitive oriented
graph $T$ with vertices $\{u_i\mid i\in [t]\}$ and strong quasi-transitive digraphs $H_1, H_2, \dots, H_t$ such that $D = T[H_1, H_2, \dots, H_t]$, where $H_i$ is substituted for $u_i$, $i\in [t]$.
\item[(b)] If $D$ is strong, then there exists a strong semicomplete
digraph $S$ with vertices $\{v_j\mid j\in [s]\}$ and quasi-transitive digraphs $Q_1, Q_2, \dots, Q_s$ such that $Q_j$ is either a vertex or is non-strong and $D = S[Q_1, Q_2, \dots, Q_s]$, where $Q_j$ is substituted for $v_j$, $j\in [s]$.
\end{description}
\end{thm}

A digraph $D$ is {\em locally semicomplete} if $N^-(x)$ and $N^+(x)$ induce semicomplete digraphs for every vertex $x$ of $D$. A {\em locally tournament digraph} or {\em local tournament} is a locally semicomplete digraph without a 2-cycle. A digraph on $n$ vertices is {\em round} if we can label its vertices $u_1, \dots, u_n$ such that for each $i$, we have $N^+(u_i)=\{u_j\mid i+1\leq j\leq i+d^+(u_i)\}$ and $N^-(u_i)=\{u_j\mid i-d^-(u_i)\leq j\leq i-1\}$ (all subscripts are taken modulo $n$). Note that every round digraph is locally semicomplete. A locally semicomplete digraph $D$ is {\em round decomposable} if there exists a round local tournament $R$ on $r\geq 2$ vertices such that $D=R[S_1, \dots, S_r]$, where each $S_i$ is a strong semicomplete digraph. We call $R[S_1, \dots, S_r]$ a {\em round decomposition} of $D$. Clearly, round decomposable locally semicomplete digraphs are totally $\Phi_2$-decomposable, where $\Phi_2$ is the set of all semicomplete digraphs and round digraphs. By an {\em evil locally semicomplete digraph} we mean a locally semicomplete digraph that is not semicomplete and not round decomposable. It is worth mentioning that Bang-Jensen, Guo, Gutin and Volkmann used the notion of round decomposition to give a full classification of locally semicomplete digraphs \cite{Bang-Jensen-Guo-Gutin-Volkmann}. 

\begin{thm}\label{intro03}\cite{Bang-Jensen-Guo-Gutin-Volkmann}
Let $D$ be a locally semicomplete digraph. Then exactly one of the
following possibilities holds. Furthermore, given a locally semicomplete digraph $D$, we can decide in polynomial time which of the possibilities holds for $D$:
\begin{description}
\item[(a)] $D$ is round decomposable with a unique round decomposition
$D=R[S_1,\\ \dots, S_r]$, where $R$ is a round local tournament on $r\geq 2$ vertices and each $S_i$ is a strong semicomplete digraph for $i\in [r]$.
\item[(b)] $D$ is evil.
\item[(c)] $D$ is a semicomplete digraph that is not round decomposable.
\end{description}
\end{thm}


The {\em lexicographic product} $G\circ H$ of two digraphs $G$ and
$H$ is a digraph with vertex set $V(G\circ H)=V(G)\times V(H)=\{(x,
x')\mid x\in V(G), x'\in V(H)\}$ and arc set $A(G\circ
H)=\{(x,x')(y,y')\mid xy\in A(G),~or~x=y~and~x'y'\in
A(H)\}$ (e.g., \cite{Hammack}). By definition, we directly have $G\circ H=G[H,\dots, H]$.

Recall that digraph compositions generalize several families of digraphs, including extended semicomplete digraphs, quasi-transitive digraphs, round decomposable locally semicomplete digraphs and lexicographic product digraphs. In particular, semicomplete compositions generalize strong quasi-transitive digraphs. To see that strong compositions form a significant generalization of strong quasi-transitive digraphs, observe that the Hamiltonian cycle problem is polynomial-time solvable for quasi-transitive digraphs \cite{Gutin4}, but NP-complete for strong semicomplete compositions (see, e.g., \cite{Bang-Jensen-Gutin-Yeo}).
While digraph composition has been used since 1990s to study locally semicomplete digraphs, quasi-transitive digraphs and their generalizations, see, e.g., \cite{Bang-Jensen-Guo-Gutin-Volkmann, Bang-Jensen2, Bang-Jensen-Huang1}, the study of digraph compositions in their own right was initiated only recently by Sun, Gutin and Ai in \cite{Sun-Gutin-Ai}.

In this paper we overview known results and get several new results on digraph compositions. The mainbody of this survey is divided into six sections: connectivity and linkage, kings and kernels, paths, cycles, acyclic spanning subdigraphs, strong spanning subdigraphs. Some conjectures and open problems are also posed for further study in this survey.

\subsection{Further notation and terminology}

Let $D$ be a digraph. If there is an arc from a vertex $x$ to a vertex $y$ in $D$, then we say that $x$ {\em dominates} $y$ and denote it by $x\rightarrow y$. If there is no arc from $x$ to $y$ we shall use the notation $x\not\rightarrow y$. If $A$ and $B$ are two subdigraphs of $D$ and every vertex of $A$ dominates each vertex of $B$, then we say that $A$ {\em dominates} $B$ and denote it by $A\rightarrow B$. We shall use $A\Rightarrow B$ to denote that $A$ dominates $B$ and there is no arcs from $B$ to $A$. A digraph $D$ is {\em triangular} with a partition $\{V_0, V_1, V_2\}$, if $V(D)$ can be partitioned into three disjoint sets $V_0, V_1, V_2$ with $V_0\Rightarrow V_1\Rightarrow V_2\Rightarrow V_0$. We use $N^-(x)$ (resp. $N^+(x)$) to denote the set of all in-neighbours (resp. out-neighbours) of a vertex $x$ in a digraph $D$.

A {\em $k$-king} in a digraph $D$ is a vertex which can reach every other vertex by a directed path of length at most $k$. A $k$-king $x$ is {\em strict} if there exists a vertex $y$ such that $d_D(x, y)=k$, where $d_D(x, y)$ denotes the length of a shortest path from $x$ to $y$ in $D$. A 2-king is called a {\em king}, and a {\em non-king} is a vertex which is not a 3-king. We call a vertex $x$ in a digraph a {\em sink} (resp. {\em source}) if $d^+(x)=0$ (resp. $d^-(x)=0$).

A set $S$ of vertices in a digraph $D$ is {\em independent} if $D[S]$ has no arcs. A subset $K$ is {\em $k$-independent} if for every pair of vertices $x,y \in K$, we have $d(x, y), d(y, x)\geq k$; it is called {\em $\ell$-absorbent} if for every $x\in V(D)\setminus K$ there exists $y\in K$ such that $d(x, y)\leq \ell$. A set $K\subseteq V(D)$ is call a {\em kernel} of $D$ if it is an independent set such that every vertex in $V(D)\setminus K$ dominates some vertex in $K$. A set $Q\subseteq V(D)$ is call a {\em quasi-kernel} of $D$ if it is an independent set such that, for every vertex $x\in V(D)\setminus Q$, there exist $y\in V(D)\setminus Q$ and $z\in Q$ such that either $xz\in A(D)$ or $xy, yz\in A(D)$.
A {\em $k$-kernel} of $D$ is a $k$-independent and $(k-1)$-absorbent subset of $V(D)$. Clearly, a 2-kernel is exactly a kernel. The problem {\sc $k$-Kernel} is determining whether a given digraph has a $k$-kernel.

The {\em complementary graph} of a graph $G$ is denoted by $\overline{G}$. A digraph $D$ is said to be {\em strong}, if for every pair of vertices $x$ and $y$, there is a path from $x$ to $y$ in $D$ and vice versa. A {\em strong component} of a digraph $D$ is a maximal induced subdigraph of $D$ which is strong. An {\em initial} (resp. {\em terminal}) {\em strong component} is one which has no arcs entering (resp. leaving) it. We call a set $S\subseteq V(D)$ a {\em separator} of $D$ if $D-S$ is not strong. A separator $S$ is {\em minimal} if no proper subset of $S$ is a separator. A digraph is {\em $k$-strong-connected} if every separator has at least $k$ vertices. 


An {\em out-tree} (resp. {\em in-tree}) {\em rooted at a vertex $r$} is an orientation of a tree such that the in-degree (resp. out-degree) of every vertex but $r$ equals one. An {\em out-branching} $B^+_r$ (resp. {\em in-branching} $B^-_r$) in a digraph $D$ is a spanning subdigraph of $D$ which is out-tree (resp. in-tree).

A digraph $D$ is {\em pancyclic} if $D$ contains a cycle of length $k$ for each $3\leq k\leq n$, and is {\em vertex-pancyclic} if every vertex of $D$ is contained in a cycle of length $k$ for each $3\leq k\leq n$. A cycle (path) of a digraph $D$ is {\em Hamiltonian} if it contains all the vertices of $D$. A digraph is {\em Hamiltonian} if it has a Hamiltonian cycle. A {\em $k$-path-cycle subdigraph} $\mathcal{F}$ of a digraph $D$ is a collection of $k$ paths $P_1, \dots, P_k$ and $t$ cycles $C_1, \dots, C_t$ such that all of $P_1, \dots, P_k, C_1, \dots, C_t$ are pairwise disjoint (possibly, $k=0$ or $t=0$). We will denote $\mathcal{F}$ by $\mathcal{F}=P_1\cup \dots \cup P_k\cup C_1\cup \dots \cup C_t$. Furthermore, a {\em $k$-path-cycle factor} is a spanning $k$-path-cycle subdigraph. If $t=0$, $\mathcal{F}$ is a {\em $k$-path subdigraph} and it is a {\em $k$-path factor} (or just a {\em path factor}) if it is spanning. If $k=0$, we say that $\mathcal{F}$ is a {\em $t$-cycle subdigraph} (or just a {\em cycle subdigraph}) and it is a {\em $t$-cycle factor} (or just a {\em cycle factor}) if it is spanning. The {\em path covering number} of a digraph $D$, denoted $pc(D)$, is the smallest $k$ for which $D$ has a {\em $k$-path factor}. Note that a Hamiltonian path is a 1-path factor. For every digraph $D$ with at least one cycle and every non-negative integer $j$, we define $\eta_j(D)=\min \{i\mid D$~has~a~$i$-path-$j$-cycle~factor$\}$.

Let $\mathcal{F}=C_1\cup \dots \cup C_t$ be a $t$-cycle factor of a digraph $D$. We say $\mathcal{F}$ is {\em reducible} if there exists a $t'$-cycle factor $\mathcal{F'}=C_1\cup \dots \cup C_{t'}$ of $D$ such that each of the following holds:\\
$(a)$~$t'<t$;\\
$(b)$~for every $i\in [t]$ there is a $j\in [t']$ such that $V(C_i)\subseteq V(C'_j)$.\\
Such a $\mathcal{F'}$ is called a {\em reduction} of $\mathcal{F}$. If no reduction of $\mathcal{F}$ exists, then $\mathcal{F}$ is said to be {\em irreducible}. Clearly, every minimum cycle factor is irreducible.

A digraph is {\em locally in-semicomplete} (resp. {\em locally out-semicomplete}) if $N^-(x)$ (resp. $N^+(x)$) induces a semicomplete digraph for every vertex $x$ of $D$. A digraph is {\em semicomplete multipartite} if it is obtained from a complete multipartite graph by replacing every edge by an arc or a pair of opposite arcs. A {\em multipartite tournament} is a semicomplete multipartite digraph without a 2-cycle. A digraph is called {\em symmetric} if for every arc $xy$ there is an opposite arc $yx$. Let $\overline{K}_p$ stand for the digraph of order $p$ with no arcs, and $\overrightarrow{C}_k$ and $\overrightarrow{P}_k$ denote the cycle and path with $k$ vertices, respectively.

\section{Connectivity and linkages}

\subsection{Connectivity}\label{sec:connectivity}

Let $\mathcal{T}_1$ be the set of all semicomplete digraphs $T$ satisfying the following: there exists a vertex $u$ such that $uv, vu\in A(T)$ for any $v\in V(T)\setminus \{u\}$. The following result concerns the structure of a minimal separator in a strong semicomplete composition $Q=T[H_1, \dots, H_t]$ when $T\not \in \mathcal{T}_1$.

\begin{thm}\label{connectivity1}\cite{Sun2}
Let $Q=T[H_1, \dots, H_t]$ be a strong semicomplete composition with $T\not \in \mathcal{T}_1$. 
Then every minimal separator $S$ of $Q$ induces in $\overline{U(Q)}$ a subgraph which consists of some connected components of $\overline{U(Q)}$. Moreover, each vertex $s\in S$ is adjacent to every vertex of $Q-S$.
\end{thm}

Note that the above result may not hold when $T\in \mathcal{T}_1$. Consider the following example \cite{Sun2}: let $T$ be a semicomplete digraph such that for each $j\in [t]$, $u_1u_j, u_ju_1\in A(T)$ and $H_i$ is a cycle of length four: $u_{i,1}, u_{i,2}, u_{i,3}, u_{i,4}, u_{i,1}$. Observe that $S=\{u_{1,1}\}\cup (\bigcup_{i=2}^t{V(H_i)})$ is a minimal separator of $Q$, and $\overline{U(Q)}$ consists of $2t$ components (each component is an edge). Clearly, the edge $\{u_{1,1}, u_{1,3}\}$ is one such component, however, $u_{1,1}\in S$ and $u_{1,3}\not\in S$.

In the same paper, Sun and Jin continued to study the connectivity of semicomplete compositions as follows.

\begin{thm}\label{connectivity2}\cite{Sun2}
Let $Q=T[H_1, \dots, H_t]$ be a $k$-strong-connected semicomplete composition with $T\not \in \mathcal{T}_1$. If $V(H_i)$ induces a connected component of $\overline{U(Q)}$ for some $i\in [t]$, then deleting all arcs in $H_i$ results in a $k$-strong-connected semicomplete composition.
\end{thm}

Note that Theorem~\ref{connectivity2} may not hold if $T\in \mathcal{T}_1$. Consider the following example \cite{Sun2}: let $T$ be a 2-cycle, and $H_i$ be a $2k$-cycle for each $1\leq i\leq 2$, where $k\geq 3$. Clearly, $T\in \mathcal{T}_1$. Observe that each $V(H_i)$ induces a connected component in $\overline{U(Q)}$ for $1\leq i\leq 2$. Indeed, each such component consists of two vertex-disjoint cliques of size $k$ and edges between them. Let $Q'$ be the subdigraph of $Q$ by deleting all arcs in $H_1$. It can be checked that $V(H_2)$ is a separator of size $k$ in $Q'$. However, there is no separator of size $k$ in $Q'$. Hence, $Q$ is $(k+1)$-strong-connected, but $Q'$ is not $(k+1)$-strong-connected.

The following result on strong quasi-transitive digraphs was given by Bang-Jensen and Huang.

\begin{thm}\label{connectivity01}\cite{Bang-Jensen-Huang1}
Every strong quasi-transitive digraph $D$ with at least four vertices has two distinct vertices $v_1, v_2$ such that $D-v_i$ is strong for $i\in [2].$
\end{thm}

Sun and Jin extended Theorem~\ref{connectivity01} to strong semicomplete compositions.

\begin{thm}\label{connectivity3}\cite{Sun2}
Every strong semicomplete composition $Q=T[H_1, \dots, H_t]$ with at least four vertices has two distinct vertices $v_1, v_2$ such that $Q-v_i$ is strong for $i\in [2].$
\end{thm}

\subsection{Linkages}\label{sec:linkage}

Let $s_1, s_2, \dots, s_k, t_1, t_2, \dots, t_k$ be distinct vertices of a digraph $D$. A {\em $k$-linkage} from $(s_1, s_2, \dots, s_k)$ to $(t_1, t_2, \dots, t_k)$ in $D$ is a set of vertex-disjoint paths $P_1,\dots ,P_k$ such that $P_i$ is an $(x_i, y_i)$-path in $D$. A digraph $D$ is {\em $k$-linked} if it contains a $k$-linkage from $(s_1, s_2, \dots, s_k)$ to $(t_1, t_2, \dots, t_k)$ for every choice of distinct vertices $s_1, s_2, \dots, s_k, t_1, t_2, \dots, t_k$.

\begin{thm}\label{link1}\cite{Bang-Jensen-ADM1989}
Every 5-strong semicomplete digraph is 2-linked and this is best possible even for tournaments.
\end{thm}

\begin{lem}\label{link2}\cite{Bang-Jensen-DM1999}
Let $Q=\overrightarrow{C}_2[H_1, H_2]$ where $H_i$ is an arbitrary digraph on $n_i$ vertices for $i\in [2]$. If $Q$ is 4-strong, then $Q$ is 2-linked.
\end{lem}

\begin{lem}\label{link3}\cite{Bang-Jensen-DM1999}
Let $k\geq 4$ be an integer and let $T$ be a digraph on at least two vertices with the property that every $k$-strong digraph of the form $T[S_1, S_2, \dots, S_t]$ is 2-linked, where $S_i$ is a semicomplete digraph on $n_i$ vertices for each $i\in [t]$. Let $Q=T[H_1, H_2, \dots, H_t]$, where $H_i$ is an arbitrary digraph on $n_i$ vertices. If $Q$ is $k$-strong, then $Q$ is 2-linked.
\end{lem}

Bang-Jensen \cite{Bang-Jensen-DM1999} deduced that every 5-strong quasi-transitive digraph is 2-linked. We extended this result to 5-strong semicomplete compositions by Theorem~\ref{link1}, Lemmas~\ref{link2} and~\ref{link3}.

\begin{thm}\label{link4}
Every 5-strong semicomplete composition is 2-linked.
\end{thm}
\begin{pf} Let $Q=T[H_1, H_2, \dots, H_t]$ be a 5-strong semicomplete composition. If $|T|=2$, then the conclusion holds by Lemma~\ref{link2}. It remains to consider the case that $|T|\geq 3$. Observe that for semicomplete digraphs $S_1, \dots, S_t$, the digraph $Q'=T[S_1, \dots, S_t]$ is semicomplete. By Theorem~\ref{link1}, if $Q'$ is 5-strong, then $Q'$ is 2-linked. Furthermore, by Lemma~\ref{link3}, $Q$ is 2-linked. 
\end{pf}

\section{Kings and kernels}

\subsection{Kings}\label{sec:king}

The study of kings in digraphs began with Landau \cite{Landau} and the term king was introduced by Maurer \cite{Maurer}. The concept of $k$-king was first introduced in \cite{Bang-Jensen-Huang2} when the authors studied quasi-transitive digraphs and this subject has received a lot of attention lately. Many of the nice existing results for kings in quasi-transitive digraphs can be naturally generalized to $k$-quasi-transitive digraphs which are exactly quasi-transitive digraphs when $k=2$\cite{Bang-Jensen-Gutin}. 
Since digraph compositions also generalize quasi-transitive digraphs (by Theorem~\ref{intro01}), it is natural to study $k$-kings of this class of digraphs.

\subsubsection{$T$ is arbitrary}

For any integer $k\geq 2$, Sun characterized all digraph compositions with a $k$-king and digraph compositions all of whose vertices are $k$-kings.

\begin{thm}\label{king1}\cite{Sun}
Let $k\geq 2$ be an integer, and $Q=T[H_1, \dots, H_t]$ be a digraph composition. Then the following assertions hold:
\begin{description}
\item[(a)] $Q$ has a $k$-king if and only if $T$ has a $k$-king $u_{i}$ for some $i\in [t]$ and at least one of the following possibilities holds: $(i)$~$H_{i}$ has a $k$-king; $(ii)$~$|V(H_{i})|\geq 2$ and $u_{i}$ belongs to a cycle of $T$ with length at most $k$.
\item[(b)] All vertices of $Q$ are $k$-kings if and only if for each $i\in [t]$, $u_{i}$ is a $k$-king of $T$ and at least one of the following possibilities holds: $(i)$~all vertices of $H_{i}$ are $k$-kings of $H_i$;   $(ii)$~$|V(H_{i})|\geq 2$ and $u_{i}$ belongs to a cycle of $T$ with length at most $k$.
\end{description}
\end{thm}

\subsubsection{$T$ is semicomplete}

The following result can be found in literature, see e.g. Theorem~2.2.9 of \cite{Bang-Jensen-Havet}.

\begin{thm}\label{king05}
Every strong semicomplete digraph is vertex-pancyclic.
\end{thm}

\begin{thm}\label{king06}\cite{Landau}
Every tournament has a king. More precisely, every vertex with maximum out-degree is a king.
\end{thm}

Bang-Jensen and Huang obtained the following result on the existence of 3-kings in a quasi-transitive digraph.

\begin{thm}\label{Bang-Huangking}\cite{Bang-Jensen-Huang2}
Let $D$ be a quasi-transitive digraph with a 3-king. If $D$ has no source, then $D$ has at least two 3-kings.
\end{thm}

Sun gave a similar result to Theorem~\ref{Bang-Huangking} for semicomplete compositions.




\begin{thm}\label{king3}\cite{Sun}
Let $Q=T[H_1, \dots, H_t]$ be a semicomplete composition with a 3-king. If $T$ has no source, then $Q$ has at least two 3-kings. In particular, every strong semicomplete composition $Q=T[H_1, \dots, H_t]$ has a 3-king, and a vertex $v\in H_i$ is a 3-king of $Q$ if and only if $u_i$ is a 3-king of $T$. 
\end{thm}

\noindent{\bf Remark A}: According to the argument for Theorem~\ref{king3}, we actually can divide $H_1, \dots, H_t$ in a strong semicomplete composition $Q=T[H_1, \dots, H_t]$ into those for which all $v\in H_i$ are 3-kings (in this case, $u_i$ is a 3-king of $T$) and those for which no $v\in H_i$ is a 3-king.

Recall that in Theorem~\ref{Bang-Huangking}, Bang-Jensen and Huang \cite{Bang-Jensen-Huang2} proved that a quasi-transitive digraph $D$ with a 3-king contains at least two 3-kings provided that it has no source.
It is also worth noting that we cannot replace the condition ``If $T$ has no source'' in Theorem~\ref{king3} by ``If $Q$ has no source'', according to the following example \cite{Sun}: let $Q=T[H_1, \dots, H_t]$ be a semicomplete composition such that $u_i$ dominates (but is not dominated by) $u_j$ for any $1\leq i< j\leq t$, and $H_1$ is a digraph with vertex set $\{u_{1, j}\mid 1\leq j\leq 6\}$ and arc set $\{u_{1,j}u_{1,j+1}, u_{1,j+1}u_{1,j} \mid 1\leq j\leq 6\}$. It can be checked that $Q$ has no source but $u_{1,4}$ is the unique 3-king of $Q$.

The following result concerns the adjacency between 3-kings and non-kings.

\begin{thm}\label{king4}\cite{Sun}
Let $Q=T[H_1, \dots, H_t]$ be a strong semicomplete composition. Then there is an arc between every 3-king and every non-king; moreover, for every non-king $u$, there exists a 3-king $v$ such that $d_Q(u, v)>3$ and $v$ dominates $u$.
\end{thm}

Bang-Jensen and Huang \cite{Bang-Jensen-Huang2} proved that Theorem~\ref{king4} holds for a quasi-transitive digraph with a 3-king, which means that the result may hold even for a non-strong quasi-transitive digraph (as long as it has a 3-king). However, for a non-strong semicomplete composition $Q$ with a 3-king, Theorem~\ref{king4} may not hold. We just use the example after Theorem~\ref{king3}: $Q$ is not strong and $u_{1,4}$ is the unique 3-king of $Q$; furthermore, there is no arc between $u_{1,4}$ and $u_{1,1}$, and of course $u_{1,4}$ does not dominate $u_{1,1}$.

For a semicomplete composition $Q$, if there exists a semicomplete composition $Q'$ which contains $Q$ as an induced subdigraph such that the set of all 3-kings of $Q'$ is precisely $V(Q)$, then we say $Q$ can be {\em established}. Similar problems for tournaments and quasi-transitive digraphs have been studied in \cite{Bang-Jensen-Huang2, Huang-Li, Reid}.
Sun gave a sufficient condition under which a semicomplete composition can be established. It suffices to study strong semicomplete compositions with a non-king, since the case that all vertices are 3-kings is trivial.

\begin{thm}\label{king5}\cite{Sun}
Let $Q=T[H_1, \dots, H_t]$ be a strong semicomplete composition. If $T$ has a strict 3-king and every 2-king of $T$ is dominated by some strict 3-king of $T$, then $Q$ can be established.
\end{thm}

For strong semicomplete multipartite digraphs, Gutin and Yeo obtained a result on the number of 4-kings.

\begin{thm}\label{king01}\cite{Gutin-Yeo}
Every strong semicomplete multipartite digraph with at least six vertices has at least five 4-kings.
\end{thm}

The following result on semicomplete bipartite digraphs was first stated by Wang and Zhang \cite{Wang-Zhang}.

\begin{thm}\label{king02}
Let $D$ be a semicomplete bipartite digraph with a unique initial strong component. If there is no 3-king in $D$, then there are at least eight 4-kings in $D$.
\end{thm}

For multipartite tournaments, Koh and Tan got the following result.

\begin{thm}\label{king03}\cite{Koh-Tan, Koh-Tan2}
If a multipartite tournament has a unique initial strong component and no 3-king, then it has at least eight 4-kings.
\end{thm}

It is worth noting that many results for strong multipartite tournaments also hold for strong semicomplete multipartite digraphs, due to the following result by Volkmann.

\begin{thm}\label{king04}\cite{Volkmann}
Every strong semicomplete $c$-partite digraph with $c \geq 3$ contains a spanning strong oriented subdigraph.
\end{thm}

By Theorems~\ref{king01}, \ref{king02}, \ref{king03} and \ref{king04}, Sun obtained the following result which deals with the minimum number of 4-kings in a strong semicomplete composition.

\begin{thm}\label{king6}\cite{Sun}
Every strong semicomplete composition $Q=T[H_1, \dots, H_t]$ with at least six vertices has at least five 4-kings. Furthermore, if $Q$ has no 3-king, then it has at least eight 4-kings.
\end{thm}

There is a generalization of the concept of kings. A {\em monochromatic king} in an arc-coloured digraph is a vertex $x$ such that for every vertex $y$, there is a monochromatic $x-y$ path. Sands, Sauer and Woodrow \cite{Sands-Sauer-Woodrow} proved that every 2-coloured tournament $T$ has a monochromatic king. Shen \cite{Shen} showed that if no subtournament of order 3 in a given arc-coloured tournament get three different colours on its arcs, then $T$ has a monochromatic king. It is natural to consider the following question of digraph compositions.

\begin{op}\label{op-king}
Give a sufficient condition under which a given arc-coloured digraph composition contains a monochromatic king.
\end{op}

\subsection{Kernels}\label{sec:kernel}

The concept of a kernel was introduced by von Neumann and Morgenstern while studying cooperative games \cite{Neumann-Morgenstern}. After that, kernels have been studied in many areas, such as list coloring, game theory and perfect graphs \cite{Boros-Gurvich}, mathematical logic \cite{Bezem-Grabmayer-Walicki} and complexity theory \cite{Walicki-Dyrkolbotn}.





Chv\'{a}tal and Lov\'{a}sz obtained the following result on the existence of a quasi-kernel in digraphs.

\begin{thm}\label{thmkernel02}\cite{Chvatal-Lovasz}
Every digraph contains a quasi-kernel.
\end{thm}

Gutin et al. \cite{Gutin-Koh-Tay-Yeo} characterized digraphs with exactly one or two quasi-kernels. In particular, if a digraph has precisely two quasi-kernels then these two quasi-kernels are actually disjoint. This raises the question of which digraphs contains a pair of disjoint quasi-kernels.
Heard and Huang \cite{Heard-Huang} gave sufficient conditions to guarantee the existence of a pair of disjoint quasi-kernels in several classes of digraphs, including semicomplete digraphs, semicomplete multipartite digraphs, quasi-transitive digraphs and locally semicomplete digraphs.

\begin{thm}\cite{Heard-Huang}\label{thmkernel01}
Every semicomplete digraph $D$ with no sink contains two vertices
$x, y$ such that $\{x\}$ and $\{y\}$ are both quasi-kernels of $D$.
\end{thm}

By Theorems~\ref{thmkernel02} and~\ref{thmkernel01}, Sun obtained the following result on the existence of disjoint quasi-kernels in semicomplete compositions.

\begin{thm}\label{thmkernel2}\cite{Sun}
Let $Q=T[H_1, \dots, H_t]$ be a semicomplete composition such that $T$ has no sink. Then $Q$ contains a pair of disjoint quasi-kernels. In particular, every strong semicomplete composition contains a pair of disjoint quasi-kernels.
\end{thm}

The following result provides a necessary condition for a semicomplete composition to have a $k$-kernel when $k\geq 3$.

\begin{lem}\label{lemkernel1}\cite{Sun}
Let $Q=T[H_1, \dots, H_t]$ be a semicomplete composition. If $Q$ contains a $k$-kernel with $k\geq 3$, then there is a vertex $v\in V(Q)$ such that $\{v\}$ is a $(k-1)$-absorbent set of $Q-(V(H_i)\setminus \{v\})$ where $v\in V(H_i)$ for some $i\in [t]$.
\end{lem}

Furthermore, when $k\geq 4$, we can characterize strong semicomplete compositions with a $k$-kernel.

\begin{lem}\label{lemkernel4}\cite{Sun}
Let $Q=T[H_1, \dots, H_t]$ be a strong semicomplete composition and $k\geq 4$. Then $Q$ has a $k$-kernel if and only if there is a vertex $v\in V(Q)$ such that $\{v\}$ is a $(k-1)$-absorbent set of $Q-(V(H_i)\setminus \{v\})$ where $v\in V(H_i)$ for some $i\in [t]$.
\end{lem}

\begin{lem}\label{lemkernel2}\cite{Sun}
Let $Q=T[H_1, \dots, H_t]$ be a digraph composition and $v\in V(H_i)$ for some $i\in [t]$. Then $\{v\}$ is a $k$-absorbent set of $Q-(V(H_i)\setminus \{v\})$ if and only if $\{u_i\}$ is a $k$-absorbent set of $T$.
\end{lem}

By Lemmas~\ref{lemkernel4} and~\ref{lemkernel2}, Sun proved the complexity of {\sc $k$-Kernel} restricted to strong semicomplete compositions for any integer $k\geq 2$.


\begin{thm}\label{thmkernel3}\cite{Sun}
Let $k\geq 2$ be an integer. The problem {\sc $k$-Kernel} restricted to strong semicomplete compositions is NP-complete when $k\in \{2,3\}$, and is polynomial time solvable when $k\geq 4$.
\end{thm}

\section{Paths}\label{sec:path}

\subsection{Hamiltonian paths}

Gutin proved the following result on Hamiltonian paths and cycles of extended semicomplete digraphs.

\begin{thm}\label{path01}\cite{Bang-Jensen-Gutin-Huang}
An extended semicomplete digraph has a Hamiltonian path (resp. cycle) if and only if it has a 1-path-cycle factor (resp. it is strong and has a cycle factor).
\end{thm}

By the above theorem, Sun and Jin characterized those semicomplete compositions which have a Hamiltonian path.

\begin{thm}\label{path1}\cite{Sun2}
Let $Q=T[H_1, \dots, H_t]$ be a semicomplete composition. Then
$Q$ has a Hamiltonian path if and only if it has a 1-path-cycle
factor $\mathcal{F}=P\cup C_1\cup \dots \cup C_k~(k\geq 0)$ such that neither $V(P)$ nor $V(C_i)$ is completely contained in a connected component of $\overline{U(Q)}$.
\end{thm}

\subsection{Path subdigraphs}

The following structural characterization of longest cycles in extended semicomplete digraphs was extensively used.

\begin{thm}\label{path02}\cite{Bang-Jensen-Huang-Yeo}
Let $Q=T[\overline{K}_{n_1}, \dots, \overline{K}_{n_t}]$ be a strong extended semicomplete digraph. For
$i\in [t]$, let $m_i$ denote the maximum number of vertices from $\overline{K}_{n_i}$ which can be covered by a cycle subdigraph of $Q$. Then every longest cycle of $Q$ contains precisely $m_i$ vertices from $\overline{K}_{n_i}$ for each $i\in [t]$.
\end{thm}

Bang-Jensen et al. proved a similar result for $k$-path subdigraphs in an extended semicomplete digraph.

\begin{thm}\label{path03}\cite{Bang-Jensen-Nielsen-Yeo}
Let $Q=T[\overline{K}_{n_1}, \dots, \overline{K}_{n_t}]$ be an extended semicomplete digraph and let $\ell_{i,k}$ denote the maximum number of vertices of $\overline{K}_{n_i}$ that can be covered by a $k$-path subdigraph in $Q$. Then every maximum $k$-path subdigraph in $Q$ covers exactly $\ell_{i,k}$ vertices of $\overline{K}_{n_i}$ for each $i\in [t]$.
\end{thm}

By Theorem~\ref{path03}, Bang-Jensen et al. got the following result on a strong semicomplete composition when each $H_i$ is either a single vertex or a non-strong quasi-transitive digraph.

\begin{thm}\label{path04}\cite{Bang-Jensen-Nielsen-Yeo}
Let $Q=T[H_1, \dots, H_t]$ be a strong semicomplete composition with order $n$, where each $H_i$ is either a single vertex or a non-strong quasi-transitive digraph. For every $k\in [n]$ and $i\in [t]$, there exists an integer $n_{i,k}$ such that every maximum $k$-path subdigraph $\mathcal{F}$ of $Q$ satisfies $|V(H_i) \cap V(\mathcal{F})|=n_{i,k}$ and no $k$-path subdigraph of $Q$ contains more than $n_{i,k}$ vertices of $H_i$.
\end{thm}

Sun and Jin extended Theorem~\ref{path04} to all strong semicomplete compositions.

\begin{thm}\label{path2}\cite{Sun2}
Let $Q=T[H_1, \dots, H_t]$ be a strong semicomplete composition with order $n$.
For every $k\in [n]$ and $i\in [t]$, there exists an integer $n_{i,k}$ such that every maximum $k$-path subdigraph $\mathcal{F}$ of $Q$ satisfies $|V(H_i) \cap V(\mathcal{F})|=n_{i,k}$ and no $k$-path subdigraph of $Q$ contains more than $n_{i,k}$ vertices of $H_i$.
\end{thm}

The argument for Theorem~\ref{path2} also means that Theorem~\ref{path03} holds for a strong semicomplete composition. Note that the number $\ell_{i, k}$ in Theorem~\ref{path3} is exactly $n_{i, k}$ of Theorem~\ref{path2}.

\begin{thm}\label{path3}\cite{Sun2}
Let $Q=T[H_1, \dots, H_t]$ be a strong semicomplete composition and let $\ell_{i,k}$ denote the maximum number of vertices of $H_i$ that can be covered by a $k$-path subdigraph in $Q$. Then every maximum $k$-path subdigraph in $Q$ covers exactly $\ell_{i,k}$ vertices of $H_i$ for each $i\in [t]$.
\end{thm}

The following nice property holds for several classes of digraphs (such as semicomplete multipartite digraphs \cite{Gutin3}): For a digraph $D$ with a 1-path-cycle subdigraph $\mathcal{F}=P_1\cup C_1\cup \dots \cup C_t$, there exists a path $P$ of $D$ such that $V(P)=V(P_1)\cup(\bigcup_{j=1}^t{V(C_j)})$. However, this property does not hold for a general semicomplete composition, according to the following example:
Let $Q=T[H_1, H_2]$ be a strong semicomplete composition such that $H_1$ is a path and $H_2$ consists of a set of disjoint cycles $C_1, \dots, C_t$, where $|V(H_1)|\leq t-2$. Clearly, $\mathcal{F}=H_1\cup C_1\cup \dots \cup C_t$ is a 1-path-cycle subdigraph but there is no path $P$ such that $V(P)=V(H_1)\cup(\bigcup_{j=1}^t{V(C_j)})$. Hence, one may try to study the following question.

\begin{op}\label{op4}
Give a sufficient condition under which the above property holds for a semicomplete composition.
\end{op}

\subsection{Path partition conjecture}

Let $\ell(D)$ denote the order of a longest directed path in a digraph $D$. The Gallai-Roy-Vitaver Theorem \cite{Gallai, Roy, Vitaver} states that the chromatic number of the underlying graph of a digraph $D$ is at most $\ell(D)$. Laborde et al. \cite{Laborde} posed a conjecture which extends this theorem in a natural way: Every digraph $D$ contains an independent set $X$ such that $\ell(D-X)< \ell(D)$. This conjecture is a  particular instance of what is called the {\em Path Partition Conjecture} which states the following:

\begin{conj}\label{conj01}\cite{Laborde}
For every digraph $D$ and any pair of positive integers $\ell_1$ and $\ell_2$ with $\ell(D)=\ell_1+\ell_2$, there exists a partition
of $D$ into two subdigraphs $D_1$ and $D_2$, such that $\ell(D_i) \leq \ell_i$  for $i \in [2]$.
\end{conj}

The following more stronger version of Conjecture~\ref{conj01} is stated in \cite{Bondy1}.

\begin{conj}\label{conj02}\cite{Bondy1}
For every digraph $D$ and any pair of positive integers $\ell_1$ and $\ell_2$ with $\ell(D)=\ell_1+\ell_2$, there exists a partition of $D$ into two subdigraphs $D_1$ and $D_2$, such that $\ell(D_i)= \ell_i$  for $i \in [2]$.
\end{conj}

In \cite{Bang-Jensen-Nielsen-Yeo}, Conjecture~\ref{conj01} is proved for several generalizations of tournaments, including quasi-transitive, extended semicomplete and locally in-semicomplete digraphs. Arroyo and Galeana-S\'anchez \cite{Arroyo} continued the research of Path Partition Conjecture for some generalizations of tournaments. 

The following conjecture is equivalent to Conjecture~\ref{conj01}.

\begin{conj}\label{conj03}\cite{Arroyo}
For every digraph $D$ and a positive integer $q< \ell(D)$, there exists a partition $(A, B)$ of $V(D)$ such that $\ell(D[A])\leq q$ and $\ell(D[B])\leq \ell(D)-q$.
\end{conj}

Arroyo and Galeana-S\'anchez proved the following results on Conjecture~\ref{conj03} for digraph compositions.

\begin{thm}\label{path4}\cite{Arroyo}
Let $Q=T[H_1, \dots, H_t]$ be a digraph composition and $q$ be a positive integer less than $\ell(D)$, where $T$ is an acyclic digraph. Suppose that for each $H_i$ and for each positive integer $q_i$ with $q_i<\ell(H_i)$, there exists a partition $(A_i, B_i)$ of $V(H_i)$ such that $\ell(H_i[A_i])\leq q_i$ and $\ell(H_i[B_i])\leq \ell(H_i)-q_i$. Then there exists a partition $(A, B)$ of $V(Q)$ such that $\ell(Q[A])\leq q$ and $\ell(Q[B])\leq \ell(Q)-q$.
\end{thm}

\begin{thm}\label{path5}\cite{Arroyo}
Let $H_1=\overline{K}_{n_1}, H_t=\overline{K}_{n_t}$, and $H_i~(2\leq i\leq t-1)$ be digraphs satisfying that for each positive integer $q_i$ with $q_i<\ell(H_i)$, there exists a partition $(A_i, B_i)$ of $V(H_i)$ with the property that $\ell(H_i[A_i])\leq q_i$ and $\ell(H_i[B_i])\leq \ell(H_i)-q_i$. Consider the composition $Q=\overrightarrow{C}_t[H_1, \dots, H_t]$ and a positive integer $q< \ell(D)$. Then there exists a partition $(A, B)$ of $V(Q)$ such that $\ell(Q[A])\leq q$ and $\ell(Q[B])\leq \ell(Q)-q$.
\end{thm}

\section{Cycles}\label{sec:cycle}

\subsection{Hamiltonian cycles}

With a similar argument to that of Theorem~\ref{path1}, Sun and Jin provided the first characterization of Hamiltonian semicomplete compositions which extends a characterization of Hamiltonian quasi-transitive digraphs given by Bang-Jensen and Huang \cite{Bang-Jensen-Huang1}.



\begin{thm}\label{cycle1}\cite{Sun2}
Let $Q=T[H_1, \dots, H_t]$ be a semicomplete composition. Then $Q$ has a Hamiltonian cycle if and only if it is strong and contains a cycle factor
$\mathcal{F}=C_1\cup \dots \cup C_k$, such that no $V(C_i)$ is completely contained in a connected component of $\overline{U(Q)}$.
\end{thm}

Gutin obtained the following result on long cycles in strong extended semicomplete digraphs.

\begin{thm}\label{Gutin-longcycle}\cite{Gutin2}
Let $D$ be a strong extended semicomplete digraph and let $\mathcal{F}$ be a cycle subdigraph of $D$. Then $D$ has a cycle $C$ which contains all vertices of $\mathcal{F}$. 
In particular, if $V(\mathcal{F})$ is maximum, then $V(C)=V(\mathcal{F})$ and $C$ is a longest cycle of $D$.
\end{thm}

By Theorem~\ref{Gutin-longcycle}, Sun and Jin extended a characterization of Hamiltonian quasi-transitive digraphs which is given by Gutin \cite{Gutin4} to strong semicomplete compositions, and provide the second characterization of Hamiltonian semicomplete compositions.

\begin{thm}\label{cycle3}\cite{Sun2}
Let $Q=T[H_1, \dots, H_t]$ be a strong semicomplete composition. Then $Q$ has a Hamiltonian cycle if and only if the extended semicomplete digraph $Q'=T[\overline{K}_{n_1}, \dots, \overline{K}_{n_t}]$ has a cycle subdigraph which covers at least $pc(H_i)$ vertices of $\overline{K}_{n_i}$ for every $i \in [t]$.
\end{thm}

\subsection{Pancyclicity}

Gutin characterized pancyclic and vertex-pancyclic extended semicomplete digraphs.
\begin{thm}\label{cycle01}\cite{Gutin}
Let $D$ be a Hamiltonian extended semicomplete digraph on $n\geq 4$ vertices such that $\overline{U(D)}$ has $k\geq 3$ connected components. Then the following assertions hold:
\begin{description}
\item[(a)] $D$ is pancyclic if and only if $D$ is not triangular with a
partition $\{V_0, V_1, V_2\}$, two of which induce independent sets, such that either $|V_0|=|V_1|=|V_2|$ or no $D[V_i]$ $(i= 0, 1, 2)$ contains a path of length 2.
\item[(b)] $D$ is vertex-pancyclic if and only if it is pancyclic and either $k>3$
or $k = 3$ and $D$ contains two cycles $C, C'$ of length 2 such that $C \cup C'$
has vertices in each of the three connected components of $\overline{U(D)}$.
\end{description}
\end{thm}

The next two lemmas by Bang-Jensen and Huang \cite{Bang-Jensen-Huang1} concern cycles in triangular digraphs.

\begin{lem}\label{cycle02}\cite{Bang-Jensen-Huang1}
Suppose that $D$ is a Hamiltonian triangular digraph with a partition $\{V_0, V_1, V_2\}$. If $D[V_1]$ contains an arc $xy$ and $D[V_2]$ contains an arc $uv$, then every vertex of $V_0\cup \{x, y, u, v\}$ is on cycles of lengths $3, 4, \dots, |V(D)|$.
\end{lem}

\begin{lem}\label{cycle03}\cite{Bang-Jensen-Huang1}
Suppose that $D$ is a triangular digraph with a partition $\{V_0, V_1, V_2\}$ and has a Hamiltonian cycle $C$. If $D[V_0]$ contains an arc of $C$ and a path $P$ of length 2, then every vertex of $V_1\cup V_2\cup V(P)$ is on cycles of lengths $3, 4, \dots, |V(D)|$.
\end{lem}

By Theorem~\ref{cycle01}, Lemmas~\ref{cycle02} and \ref{cycle03}, Sun and Jin got the following result on the pancyclicity of semicomplete compositions, and this extends a similar result for quasi-transitive digraphs by Bang-Jensen and Huang \cite{Bang-Jensen-Huang1}.

\begin{thm}\label{cycle2}\cite{Sun2}
Let $Q=T[H_1, \dots, H_t]$ be a Hamiltonian semicomplete composition on $n\geq 4$ vertices. Then $Q$ is pancyclic if and only if it is not triangular with a partition
$\{V_0, V_1, V_2\}$, two of which induce independent sets, such that either $|V_0|=|V_1|=|V_2|$, or no $Q[V_i]~(i=0, 1, 2)$ contains a path of length 2.
\end{thm}

Theorem~\ref{cycle2} completely characterizes pancyclic Hamiltonian semicomplete compositions, but we still have no idea for the characterization of vertex-pancyclic Hamiltonian semicomplete compositions.

\begin{op}\label{op2}
Characterize Hamiltonian semicomplete compositions which are vertex-pancyclic.
\end{op}

Recall that Theorem~\ref{Gutin-longcycle} states that a strong extended semicomplete digraph $D$ has a cycle $C$ which contains all vertices of $\mathcal{F}$, where $\mathcal{F}$ is a cycle subdigraph of $D$. In particular, if $V(\mathcal{F})$ is maximum, then $V(C)=V(\mathcal{F})$ and $C$ is a longest cycle of $D$. Gutin also got the following result on strong semicomplete bipartite digraphs. 

\begin{thm}\label{Gutin-1}\cite{Gutin2}
Let $D$ be a strong semicomplete bipartite digraph. The length of a longest cycle in $D$ is equal to the number of vertices in a cycle subdigraph of $D$ of maximum order. 
\end{thm}


However, it seems that Theorems~\ref{Gutin-longcycle} and~\ref{Gutin-1} do not hold for strong semicomplete compositions, according to the following example \cite{Sun2}: Let $Q=T[H_1, H_2]$ be a strong semicomplete composition such that $|V(H_2)|=|V(H_1)|+2$, $H_1$ has no arcs, $H_2$ contains a 2-cycle and has no other arcs. Since $Q$ is strong and $t=2$, $T$ is a 2-cycle. Observe that the order of a maximum cycle subdigraph of $Q$ is $|V(Q)|$, but the length of the longest cycle is $|V(Q)|-1$. Therefore, it is interesting to study the following question.

\begin{op}\label{op1}
Let $Q=T[H_1, \dots, H_t]$ be a strong semicomplete composition.
Give a sufficient condition under which the length of a longest cycle in $Q$ is equal to the number of vertices in a cycle subdigraph of $Q$ of maximum order (or, $Q$ has a cycle $C$ which contains all vertices of a given cycle subdigraph $\mathcal{F}$ of $Q$).
\end{op}

\subsection{Cycle factors}

For a digraph composition $Q=T[H_1, H_2, \dots, H_t]$, those cycles of a cycle factor $\mathcal{F}$ that are contained in a $H_i$ are
called {\em small} cycles and all other cycles of $\mathcal{F}$ are called {\em large} cycles. Using Theorem~\ref{path02}, We prove the following result:

\begin{lem}\label{cyclefactor1}
Let $\mathcal{F}$ be an irreducible cycle factor in a strong semicomplete composition $Q=T[H_1, H_2, \dots, H_t]$. Then $\mathcal{F}$ has exactly one large cycle. In particular, every minimum cycle factor in a strong semicomplete composition contains exactly one large cycle.
\end{lem}
\begin{pf}
Let $\mathcal{F}=C_1\cup \dots \cup C_q$ be an irreducible $q$-cycle factor of $Q$. Without loss of generality, assume that $C_1, \dots, C_p$ are all large cycles, where $0\leq p\leq q$. Let $Q^*$ denote the strong extended semicomplete digraph which is obtained from $Q$ by contracting every maximal subpath of $C_j$ inside every $H_i$ and then deleting all remaining arcs in $H_i$, where $i\in [t], j\in [q]$. $Q^*$ must be of the form $T[\overline{K}_{n'_1}, \dots, \overline{K}_{n'_t}]$, and therefore is strong, where $1\leq n'_i\leq n_i$ for each $i\in [t]$. Moreover, after the above operations, every small cycle of $\mathcal{F}$ is contracted into a vertex and the set of large cycles $\{C_i\mid i\in [p]\}$ is converted into a cycle subdigraph $\mathcal{F}^*$ of $Q^*$. Let $C^*$ be a longest cycle of $Q^*$. By Theorem~\ref{path02}, we can construct $C^*$ such that it covers at least one vertex from every $V(\overline{K}_{n'_i})$ (since now each vertex of $V(\overline{K}_{n'_i})$ belongs to a cycle of $Q^*$ by the fact that $Q^*$ is strong), every vertex of $\mathcal{F}^*$ and possibly some of the other vertices (those corresponding to small cycles in $\mathcal{F}$).

Now we can obtain a new cycle subdigraph $\mathcal{F}'$ which contains one large cycle, say $C'$, corresponding to $C^*$, and possibly some small cycles all of which belong to $\mathcal{F}$, by substituting the contracted paths (including contracted small cycles). Observe that when $p=0$ or $p>1$, $\mathcal{F}'$ has fewer cycles than $\mathcal{F}$ and the vertices of each cycle of $\mathcal{F}$ are covered by one cycle of $\mathcal{F}'$, which means that $\mathcal{F}'$ is a reduction of $\mathcal{F}$, a contradiction. Therefore, $\mathcal{F}$ has precisely one large cycle. In particular, every minimum cycle factor in a strong semicomplete composition contains exactly one large cycle.
\end{pf}

Let $\mathcal{C}$ be the set of cycle subdigraphs of $T[\overline{K}_{n_1}, \dots, \overline{K}_{n_t}]$ which is the corresponding extended semicomplete digraph of a strong semicomplete composition $Q=T[H_1, H_2, \dots, H_t]$. For each $i\in [t]$, let $$m_i(Q)=\max_{\mathcal{S}\in \mathcal{C}}\{|V(\mathcal{S})\cap V(\overline{K}_{n_i})|\}.$$
By Theorem~\ref{path02}, every longest cycle in $T[\overline{K}_{n_1}, \dots, \overline{K}_{n_t}]$ contains exactly $m_i(Q)$ vertices of $V(\overline{K}_{n_i})$.

\begin{lem}\label{cyclefactor2}
If $Q=T[H_1, H_2, \dots, H_t]$ is a strong semicomplete composition containing a cycle factor then $Q$ has a minimum cycle factor in
which the (unique) large cycle $C$ intersects $H_i$ in exactly $m_i(Q)$ paths for each $i\in [t]$. That is, by contracting
each maximal subpath of $C$ which lies inside $H_i$ (for every $i\in [t]$), we obtain a longest cycle of $T[\overline{K}_{n_1}, \dots, \overline{K}_{n_t}]$.
\end{lem}
\begin{pf}
Let $\mathcal{F}$ be a minimum cycle factor in $Q$ with the large cycle $C_0$. By Lemma~\ref{cyclefactor1}, $C_0$ is unique. For each $i\in [t]$, let $p_i$ be the number of maximal subpaths of $C_0$ inside $H_i$. Observe that $p_i\leq m_i(Q)(\leq n_i)$.

If $H_i$ contains no small cycle from $\mathcal{F}$, then $H_i$ clearly contains an $m_i(Q)$-path cover, since now we have $p_i\geq pc(H_i)$.

Otherwise, assume that $H_i$ contains $c_i(\geq 1)$ small cycles from $\mathcal{F}$, we claim that $p_i=m_i(Q)$ in this case. Suppose not, that is, $p_i< m_i(Q)$. We delete one arc in a small cycle in $Q_i$ until we either get $m_i(Q)$ paths or have used all small cycles (in this case, we also delete some arcs in the current $p_i+c_i$ paths until we have an $m_i(Q)$-path cover). In both cases, by doing this in each $H_i$ and replacing $C_0$ by a large cycle entering $H_i$ exactly $m_i(Q)$ times (note that such a large cycle can be constructed from a longest cycle in $T[\overline{K}_{n_1}, \dots, \overline{K}_{n_t}]$ which contains exactly $m_i(Q)$ vertices of $V(\overline{K}_{n_i})$, by replacing these $m_i(Q)$ vertices with the $m_i(Q)$ paths in an $m_i(Q)$-path cover), we obtain a cycle factor with fewer cycles than $\mathcal{F}$, this produces a contradiction.

Now we can construct a minimum cycle factor with the desired property from $\mathcal{F}$ and any longest cycle $C$ in $T[\overline{K}_{n_1}, \dots, \overline{K}_{n_t}]$ by the following operation: for the $m_i(Q)$ vertices of $V(\overline{K}_{n_i})\cap C$, substitute $m_i(Q)$ paths of $H_i$ and keep the small cycles of $\mathcal{F}$.    
\end{pf}

A {\em canonical} minimum cycle factor is one for which the unique large cycle intersects each $H_i$ in exactly $m_i(Q)$ paths. By Lemma~\ref{cyclefactor2}, every strong semicomplete composition with a cycle factor has a canonical minimum cycle factor.

Let $I(Q)=\{i\mid m_i(Q)<pc(H_i)\}$ and note that for every $i$, $H_i$ has the same number $c_i$ of small cycles with respect to every canonical minimum cycle factor. Observe that the number $c_i=0$ when $i\not\in I(Q)$. Furthermore, it can be checked that $c_i=\min\{j\mid \eta_j(H_i)=m_i(Q)\}$. Therefore, we have the following characterization of the number, $k_{\min}(Q)$, of cycles in a minimum cycle factor of a strong semicomplete composition.

\begin{thm}\label{cyclefactor3} Let $Q=T[H_1, H_2, \dots, H_t]$ be a strong semicomplete composition with a cycle factor. Then
$$k_{\min}(Q)=1+\sum_{i\in I(Q)}{\min\{j\mid \eta_j(H_i)=m_i(Q)\}}.$$
\end{thm}

Note that Theorem~\ref{cyclefactor3} extends a similar result for strong quasi-transitive digraphs (see Theorem~21 of \cite{Bang-Jensen-Nielsen}): Let $D$ be a strong quasi-transitive digraph with a cycle factor which has a canonical decomposition $D=T[H_1, H_2, \dots, H_t]$. Then $k_{\min}(D)=1+\sum_{i\in I(D)}{\min\{j\mid \eta_j(H_i)=m_i(D)\}}$. Bang-Jensen and Nielsen studied the problem of determining the complexity of computing $k_{\min}(D)$ and finding a minimum cycle factor of a quasi-transitive digraph $D$, and deduced that when $k_{\min}(D)$ is small, in particular when $k_{\min}(D)=\{1,2,3\}$, the above problem is polynomially solvable. However, when we consider the same problem on strong semicomplete compositions, it becomes NP-complete even when $k_{\min}(D)=1$, as in this case the problem is exactly the Hamiltonian cycle problem which was proved to be NP-complete for strong semicomplete compositions \cite{Bang-Jensen-Gutin-Yeo}.

\section{Acyclic spanning subdigraphs}\label{sec:ass}

\subsection{The existence of prescribed acyclic spanning subdigraphs}

It is well known that a tournament $T$ contains an $x-y$ Hamiltonian path if and only if there is an acyclic spanning subdigraph $R$ (not necessarily induced) such that for each vertex $z$ of $T$, $R$ contains an $x-z$ path and a $z-y$ path \cite{Thomassen2}. Bang-Jensen and Huang \cite{Bang-Jensen-Huang1} proved that if a quasi-transitive digraph has both in- and out-branchings then it always contains such an acyclic spanning subdigraph. For semicomplete compositions, we give the following sufficient condition to guarantee the existence of this type of subdigraph.

\begin{thm}\label{ass1}
Let $Q=T[H_1, \dots, H_t]$ be a semicomplete composition. Then it contains an acyclic spanning subdigraph $R$ with a source $x$ and a sink $y$ such that for each vertex $z$ of $Q$, $R$ contains an $x-z$ path and a $z-y$ path, if one of the following assertions holds:
\begin{description}
\item[(a)] $Q$ is non-strong and has both in- and out-branchings.
\item[(b)] $Q$ is strong with $|V(H_i)|\geq 2$ for each $i\in [t]$.
\end{description}
\end{thm}
\begin{pf}
We first assume that $(a)$ holds. Let $Q$ have an out-branching, say $B^+_x$, rooted at $x$ and an in-branching, $B^-_y$, rooted at $y$, then $Q$ has precisely one initial strong component, say $Q'$, and precisely one strong terminal component, say $Q''$. Furthermore, we must have $x\in V(Q')$ and $y\in V(Q'')$. By the definition of a semicomplete composition, $Q'$ (resp. $Q''$) is a subdigraph of some $H_i$ when $Q'$ (resp. $Q''$) contains vertices from only one $H_i$, or is the union of some $H_i$ when $Q'$ (resp. $Q''$) contains vertices from at least two $H_i$.

We just consider the case that $Q'$ is a subdigraph of $H_1$, and $Q''=\bigcup_{i=s}^t{H_i}$ since the other cases are similar. Since $Q'$ (resp. $Q''$) is the unique strong initial (resp. terminal) strong component, we have $V(H_1)\Rightarrow \bigcup_{i=2}^{s-1}{H_i}\Rightarrow \bigcup_{i=s}^t{H_i}$. Let $B'^+_x$ be the subdigraph of $B^+_x$ induced by $V(H_1)$, it is not hard to see that $B'^+_x$ is an out-branching of $H_1$ rooted at $x$. Since $Q''$ is strong, it has an in-branching $B'^-_y$ rooted at $y$. We now construct a subdigraph $R$ of $Q$ from $B'^+_x$ and $B'^-_y$ by adding all arcs from $V(H_1)$ to $v$ and all arcs from $v$ to $V(Q'')$ for each $v\in \bigcup_{i=2}^{s-1}{V(H_i)}$ (if $\bigcup_{i=2}^{s-1}{V(H_i)}=\emptyset$, then we just add all arcs from $V(H_1)$ to $V(Q'')$). It can be checked that $R$ is the desired acyclic spanning subdigraph.

We next assume that $(b)$ holds. By 
Theorem~\ref{king05} and the fact \cite{Sun2} that a digraph composition $Q=T[H_1, \dots, H_t]$ is strong if and only if $T$ is strong, $T$ is strong and therefore has a Hamiltonian cycle $C: u_1, u_2, \dots, u_t, u_1$.

If $t\geq 3$, we construct a digraph $R$ from $C$ as follows: For each $i\in [t]$, substitute a copy of $H_i$ for $u_i$ and then delete all arcs inside $H_i$. Let $x\in V(H_1)$ and $y\in V(H_2)$; delete all arcs from $V(H_t)$ to $x$, all arcs from $V(H_1)\setminus \{x\}$ to $V(H_2)\setminus \{y\}$ and all arcs from $y$ to $V(H_3)$. It is not hard to check that $R$ is an acyclic spanning subdigraph of $Q$ satisfying the desired properties.

For the case that $t=2$, we construct a digraph $R$ from $C$ as follows: For each $i\in [2]$, substitute a copy of $H_i$ for $u_i$ and then delete all arcs inside $H_i$. Let $x\in V(H_1)$ and $y\in V(H_2)$; delete all arcs from $V(H_1)\setminus \{x\}$ to $V(H_2)\setminus \{y\}$, from $V(H_2)$ to $x$, and all arcs from $y$ to $V(H_1)$. Observe that $R$ is the desired acyclic spanning subdigraph of $Q$.
\end{pf}

Note that we use the Hamiltonicity (Theorem~\ref{king05}) of a strong semicomplete digraph in the proof for the case that $Q$ is strong with $|V(H_i)|\geq 2$ for each $i$ in Theorem~\ref{ass1}. In fact, the proof also means that the following more general result holds.

\begin{thm}\label{ass2}
Let $Q=T[H_1, \dots, H_t]$ be a digraph composition. If $T$ is Hamiltonian and $|V(H_i)|\geq 2$ for each $i\in [t]$, then $Q$ contains an acyclic spanning subdigraph $R$ with a source $x$ and a sink $y$ such that for each vertex $z$ of $Q$, $R$ contains an $x-z$ path and a $z-y$ path.
\end{thm}

In Theorem~\ref{ass2}, the existence of prescribed acyclic spanning subdigraphs in strong semicomplete compositions with $|V(H_i)|\geq 2$ for each $i\in [t]$ were studied. It is natural to extend this to all strong semicomplete compositions. However, consider the following example: let $Q=T[H_1, H_2]$, where $T$ is a 2-cycle, $|V(H_1)|=1$, $|V(H_2)|\geq 3$ and $V(H_2)$ induces an independent set of $Q$. It can be checked that $Q$ does not contain an acyclic spanning subdigraph $R$ with a source $x$ and a sink $y$ such that for each vertex $z$ of $Q$, $R$ contains an $x-z$ path and a $z-y$ path. This means that not all strong semicomplete compositions contain such a subdigraph. Hence, it is interesting to study the following question.

\begin{op}\label{op3}
Characterize strong semicomplete compositions $Q$ which contain an acyclic spanning subdigraph $R$ with a source $x$ and a sink $y$ such that for each vertex $z$ of $Q$, $R$ contains an $x-z$ path and a $z-y$ path.
\end{op}

\subsection{Arc-disjoint in-and out-branchings}

Edmonds \cite{Edmonds} characterized digraphs with $k$ arc-disjoint out-branchings rooted at a specified vertex $r.$ Furthermore, there exists a polynomial algorithm for finding $k$ arc-disjoint out-branchings with a given root $r$ if they exist \cite{Bang-Jensen-Gutin}. However, it is NP-complete to decide whether a digraph $D$ has a pair of arc-disjoint out-branching and in-branching rooted at $r,$ which was proved by Thomassen (see \cite{Bang-Jensen}). Following \cite{Bang-Jensen-Huang2014} we will call such a pair a {\em good pair rooted at $r$}. Note that a good pair forms a strong spanning subdigraph of $D$ and thus if $D$ has a good pair, then $D$ is strong. The problem of the existence of a good pair was studied for tournaments and their generalizations, and characterizations (with proofs
implying polynomial-time algorithms for finding such a pair) were obtained in \cite{Bang-Jensen} for tournaments,  \cite{Bang-Jensen-Huang1}
for quasi-transitive digraphs and \cite{Bang-Jensen-Huang2014} for locally semicomplete digraphs. Also, Bang-Jensen and Huang \cite{Bang-Jensen-Huang1}  showed that if $r$ is adjacent to every vertex of $D$ (apart from itself) then $D$ has a good pair rooted at $r$.

In \cite{Gutin-Sun}, Gutin and Sun studied the existence of good pairs for digraph compositions. By Theorem~\ref{sss2}, they got the following lemma.

\begin{lem}\label{ass3}\cite{Gutin-Sun}
Let $Q=T[H_1,\dots ,H_t],$ where $t\ge 2.$ If $T$ has a Hamiltonian cycle and $H_1,\dots ,H_t$ are arbitrary digraphs, each with at least two vertices, then $Q$ has a good pair at any root $r$.
\end{lem}

An {\em ear decomposition} of a digraph $D$ is a sequence
$\mathcal{P}=(P_0, P_1, P_2, \cdots, P_t)$, where $P_0$ is a cycle
or a vertex and each $P_i$ is a path, or a cycle with the following
properties:\\
$(a)$~$P_i$ and $P_j$ are arc-disjoint when $i\neq j$.\\
$(b)$~For each $i=0,1,2,\cdots,t$: let $D_i$ denote the digraph with
vertices $\bigcup_{j=0}^i{V(P_j)}$ and arcs
$\bigcup_{j=0}^i{A(P_j)}$. If $P_i$ is a cycle, then it has
precisely one vertex in common with $V(D_{i-1})$. Otherwise the end
vertices of $P_i$ are distinct vertices of $V(D_{i-1})$ and no other
vertex of $P_i$ belongs to $V(D_{i-1})$.\\
$(c)$~$\bigcup_{j=0}^t{A(P_j)}=A(D)$.

The following result on ear decomposition is well-known, see e.g., \cite{Bang-Jensen-Gutin}.
\begin{thm}\label{ass01}
Let $D$ be a digraph with at least two vertices. Then $D$ is strong
if and only if it has an ear decomposition. Furthermore, if $D$ is
strong, every cycle can be used as a starting cycle $P_0$ for an ear
decomposition of $D$, and there is a linear-time algorithm to find such
an ear decomposition.
\end{thm}

By Lemma~\ref{ass3} and Theorem~\ref{ass01}, Gutin and Sun furthermore got the following result.
\begin{lem}\label{ass4}\cite{Gutin-Sun}
Let $Q=T[\overline{K}_2, \dots, \overline{K}_2],$ where $|V(T)|=t\ge 2$ and $\overline{K}_2$ is the digraph with two vertices and no arcs.
If $T$ is strong, then $Q$ has a good pair at any root $r$.
\end{lem}

By Lemmas~\ref{ass3} and~\ref{ass4}, and Theorem~\ref{ass01}, Gutin and Sun obtained a sufficient condition under which a digraph composition has a good pair at any root $r$.

\begin{thm}\label{ass5}\cite{Gutin-Sun}
Let $Q=T[H_1,\dots ,H_t],$ where $t\ge 2.$ If $T$ is strong and $H_1,\dots ,H_t$ are arbitrary digraphs, each with at least two vertices, then $Q$ has a good pair at any root $r$. Furthermore, this pair can be found in polynomial time.
\end{thm}

The condition of $n_i\ge 2$ in Theorem~\ref{ass5} cannot be relaxed. Indeed, the following characterization of quasi-transitive digraphs with a good pair provides an infinite family of strong quasi-transitive digraphs which have no good pair rooted at some vertices.

\begin{thm}\label{ass02}\cite{Bang-Jensen-Huang1}
Let $D$ be a strong quasi-transitive digraph and $r$ a vertex of $D$ such that $V(D)=\{r\}\cup N^-(r)\cup N^+(r)$. There is a polynomial-time algorithm to decide whether $D$ has a good pair at $r$.
\end{thm}

Gutin and Sun characterized semicomplete compositions with a good pair,
which generalizes the corresponding characterization in Theorem~\ref{ass02}.

\begin{thm}\label{ass6}\cite{Gutin-Sun}
A strong semicomplete composition $Q$ has a good pair rooted at $r$ if and only if $Q'=Q[\{r\}\cup N^-(r)\cup N^+(r)]$ has a good pair rooted at $r$.
\end{thm}

By Theorems \ref{ass02} and \ref{ass6},  we immediately have the following:

\begin{thm}\label{ass7}\cite{Gutin-Sun}
Given a semicomplete composition and a vertex $r$, we can decide in polynomial time whether $D$ has a good pair rooted at $r.$
\end{thm}

Recall that Theorem~\ref{ass6} generalizes a similar characterization by Bang-Jensen and Huang \cite{Bang-Jensen-Huang1} for quasi-transitive digraphs. Strong semicomplete compositions is not the only class of digraphs generalizing strong quasi-transitive digraphs. Other such classes have been studied such as $k$-quasi-transitive digraphs  and it would be interesting to see whether a characterization for the problem (or, at least non-trivial sufficient conditions) on $k$-quasi-transitive digraphs can be obtained.

\begin{op}\label{op5}\cite{Gutin-Sun}
Characterize $k$-quasi-transitive digraphs with a good pair at any root $r$.
\end{op}

As we mentioned above, Bang-Jensen and Huang \cite{Bang-Jensen-Huang2014} obtained a characterization for the problem for locally semicomplete digraphs. It would be interesting to see whether a characterization for the problem on locally in-semicomplete digraphs can be obtained.

\begin{op}\label{op6}\cite{Gutin-Sun}
Characterize locally in-semicomplete digraphs with a good pair at any root $r$.
\end{op}

It is worth mentioning a related question. An out-branching and in-branching $B^+_r$ and $B^-_r$ are called $k$-{\em distinct} if $B^+_r$ has at least $k$ arcs not present in $B^-_r.$ The problem of deciding whether a digraph $D$ has a $k$-distinct pair of out- and in-branchings is NP-complete since it generalizes the good pair problem ($k=|V(D)|-1$). Bang-Jensen and Yeo \cite{bangDAM156} asked whether the $k$-distinct problem is fixed-parameter tractable when parameterized by $k$, i.e., whether there is an $O(f(k)|V(D)|^{O(1)})$-time algorithm for solving the problem, where $f(k)$ is an arbitrary computable function in $k$ only. Gutin, Reidl and Wahlstr{\"o}m \cite{gutinJCSS95}  answered this open question in affirmative by designing an $O(2^{O(k\log^2 k)}|V(D)|^{O(1)})$-time algorithm for solving the problem.

Recently, Bang-Jensen and Wang~\cite{Bang-Jensen-Wang} introduced the concept of a {\em good $(u,v)$-pair} which is defined as a pair of arc-disjoint out-branching $B^+_u$ and in-branching $B^-_v$. They considered the existence of good $(u,v)$-pairs in compositions of strong semicomplete
digraphs and transitive digraphs and give a complete classification of semicomplete and transitive compositions with no good $(u, v)$-pair for given vertices $u, v$. Their proofs are constructive and can be converted to polynomial algorithms. 

\begin{thm}\label{thmBANG-WANG}\cite{Bang-Jensen-Wang}
There exists a polynomial algorithm which given a composition $Q=T[H_1, \dots, H_t]$, where $T$ is either transitive or semicomplete, and two vertices $u, v$, decides whether $Q$ has a good $(u, v)$-pair and outputs such a pair when it exists.    
\end{thm}

Furthermore, by Theorems~\ref{intro01} and~\ref{thmBANG-WANG}, they proved the following result which confirms a conjecture by Bang-Jensen and Gutin \cite{Bang-Jensen-Gutin1998}.

\begin{cor}\label{corBANG-WANG}\cite{Bang-Jensen-Wang}
There exists a polynomial algorithm which given a quasi-transitive digraph $D$ and two vertices $u, v$, decides whether $D$ has a good $(u, v)$-pair and outputs such a pair when it exists.    
\end{cor}

\section{Strong spanning subdigraphs}\label{sec:sss}

\subsection{Strong spanning subdigraphs without a 2-cycle}

The following result on strong semicomplete digraphs can be found in the literature, see e.g. Proposition~2.2.8 of \cite{Bang-Jensen-Havet}.

\begin{thm}\label{sss01}
Every strong semicomplete digraph on $n \geq 3$ vertices contains a strong spanning tournament.
\end{thm}

Theorem~\ref{sss01} means that every strong semicomplete digraph on $n \geq 3$ vertices contains a strong spanning semicomplete digraph without a 2-cycle. Bang-Jensen and Huang obtained a similar result for strong quasi-transitive digraphs.

\begin{thm}\label{sss02}\cite{Bang-Jensen-Huang1}
Every strong quasi-transitive digraph contains a strong spanning quasi-transitive digraph without a 2-cycle.
\end{thm}

Sun and Jin got a similar result to Theorem~\ref{sss02} for strong semicomplete compositions.

\begin{thm}\label{sss1}\cite{Sun2}
A strong semicomplete composition $Q=T[H_1, \dots, H_t]$ contains a strong spanning semicomplete composition without a 2-cycle if
and only if $t\geq 3$.
\end{thm}

\subsection{Smallest strong spanning subdigraph}

For a digraph $D$ and a natural number $k$, let $H_k(D)$ denote a digraph obtained from $D$ as follows: add two sets of $k$ new vertices $x_1, x_2, \dots, x_k$, $y_1, y_2, \dots, y_k$; add all possible arcs from $V(D)$ to $x_i$ along with all possible arcs from $y_i$ to $V(D)$, $i\in [k]$; add all arcs of the kind $x_iy_j$, $i, j \in [k]$. Clearly, $H_0(D)= D$. Let $D$ be a strong connected digraph and let $\epsilon(D)$ be the smallest $k \geq 0$ such that $H_k(D)$ is Hamiltonian. Observe that 
\[
\epsilon(D)=\left\{
   \begin{array}{ll}
     0, &\mbox {$D$ is Hamiltonian;}\\
     pc(D), &\mbox {Otherwise.}
   \end{array}
   \right.
\]

The {\em smallest strong spanning subdigraph} of a digraph $D$ is defined to be the strong spanning subdigraph of $D$ with smallest arcs. The problem of
{\sc Smallest Strong Spanning Subdigraph} (SSSS) is defined as follows: given a strong digraph $D$, find a strong spanning subdigraph $D'$ of $D$ such that $D'$ has as few arcs as possible. This problem generalizes the Hamiltonian cycle problem (and therefore is NP-hard) and has been considered in the literature, see e.g. \cite{Aho-Garey-Ullman, Bang-Jensen-Huang-Yeo, Khuller-Raghavachari-Young}.

In particular, Bang-Jensen, Huang and Yeo obtained the following result on the smallest strong spanning subdigraph of a strong quasi-transitive digraph.

\begin{thm}\label{sss03}\cite{Bang-Jensen-Huang-Yeo}
The smallest strong spanning subdigraph of a strong quasi-transitive digraph $D$ has precisely $n+\epsilon(D)$ arcs.
\end{thm}

In the argument of Theorem~\ref{sss03}, the authors used the following two lemmas.

\begin{lem}\label{BHY2003-1}\cite{Bang-Jensen-Huang-Yeo}
If $D$ is an acyclic extended semicomplete digraph, then $pc(D)= \max\{|I| \mid I~is~an~independent~set~in~D\}$. Furthermore, starting from $D$, one can obtain a path covering with $pc(D)$ paths by removing the vertices of a longest path $pc(D)$ times.
\end{lem}

\begin{lem}\label{BHY2003-2}\cite{Bang-Jensen-Huang-Yeo}
Let $D$ be a strong extended semicomplete digraph and let $C$ be a longest
cycle in $D$. Then $D-V(C)$ is acyclic.
\end{lem}

By Theorems~\ref{cycle1}, Lemmas~\ref{BHY2003-1} and~\ref{BHY2003-2}, Sun and Jin extended Theorem~\ref{sss03} to strong semicomplete compositions.


\begin{thm}\label{sss8}\cite{Sun2}
The smallest strong spanning subdigraph of a strong semicomplete composition $Q=T[H_1, \dots, H_t]$ has precisely $n+\epsilon(Q)$ arcs.
\end{thm}



\subsection{Strong arc decomposition}

Recall that it is NP-complete to decide whether a digraph $D$ has a pair of arc-disjoint out-branching and in-branching rooted at $r,$ which was proved by Thomassen (see \cite{Bang-Jensen}). In connection with this problem, Thomassen \cite{Thomassen} posed the following conjecture: There exists an integer $N$ so that every $N$-arc-strong digraph $D$ contains a pair of arc-disjoint in- and out-branchings. A digraph $D=(V,A)$ has a {\em strong arc decomposition} if $A$ has two disjoint sets $A_1$ and $A_2$ such that both
$(V,A_1)$ and $(V,A_2)$ are strong \cite{Bang-Jensen-Gutin-Yeo}. Bang-Jensen and Yeo \cite{Bang-Jensen-Yeo} generalized \footnote{Every strong digraph $D$ has an out- and in-branching rooted at any vertex of $D$.} the above conjecture as follows: There exists an integer $N$ so that every $N$-arc-strong digraph $D$
has a good decomposition. For a general digraph $D$, it is NP-complete to decide whether a digraph has a strong arc decomposition \cite{Bang-Jensen-Gutin-Yeo}. In the same paper, Bang-Jensen and Yeo characterized semicomplete digraphs with a strong arc decomposition. Bang-Jensen and Huang \cite{Bang-Jensen-Huang} extended this result by characterizing locally semicomplete digraphs with such a decomposition.
After that, Sun, Gutin and Ai \cite{Sun-Gutin-Ai} continued research on strong arc decompositions in classes of digraphs and consider digraph compositions
and products.

\begin{figure}[!h]
\begin{center}
\tikzstyle{vertexX}=[circle,draw, top color=gray!5, bottom color=gray!30, minimum size=16pt, scale=0.6, inner sep=0.5pt]
\tikzstyle{vertexY}=[circle,draw, top color=gray!5, bottom color=gray!30, minimum size=20pt, scale=0.7, inner sep=1.5pt]

\begin{tikzpicture}[scale=0.4]
 \node (v1) at (1.0,6.0) [vertexX] {$v_1$};
 \node (v2) at (7.0,6.0) [vertexX] {$v_2$};
 \node (v3) at (1.0,1.0) [vertexX] {$v_3$};
 \node (v4) at (7.0,1.0) [vertexX] {$v_4$};
\draw [->, line width=0.03cm] (v1) -- (v2);
\draw [->, line width=0.03cm] (v2) -- (v3);
\draw [->, line width=0.03cm] (v3) -- (v4);
\draw [->, line width=0.03cm] (v4) -- (v1);
\draw [->, line width=0.03cm] (v1) to [out=285, in=75] (v3);
\draw [->, line width=0.03cm] (v3) to [out=105, in=255] (v1);
\draw [->, line width=0.03cm] (v2) to [out=285, in=75] (v4);
\draw [->, line width=0.03cm] (v4) to [out=105, in=255] (v2);
\end{tikzpicture}
\caption{Digraph $S_4$}\label{S4fig}
\end{center}
\end{figure}
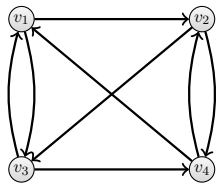

The following theorem gives sufficient conditions for a digraph composition to have a strong arc decomposition. We use $S_4$ to denote the digraph obtained from the complete digraph with four vertices by deleting a cycle of length 4, as shown in Figure~\ref{S4fig}.

\begin{thm}\label{sss2}\cite{Sun-Gutin-Ai}
Let $T$ be a digraph with vertices $u_1,\dots ,u_t$ ($t\ge 2$) and let $H_1,\dots ,H_t$ be digraphs. Let the vertex set of $H_i$ be $\{u_{i,j_i}\mid i\in [t], j_i\in [n_i]\}$ for every $i\in [t].$
Then $Q=T[H_1,\dots ,H_t]$ has a strong arc decomposition if at least one of the following conditions holds:
\begin{description}
\item[(a) ]  $T$ is a 2-arc-strong semicomplete digraph and $H_1,\dots , H_t$ are arbitrary digraphs, but $Q$ is not isomorphic to $S_4;$
\item[(b) ] $T$ has a Hamiltonian cycle and one of the following conditions holds:
\begin{itemize}
\item $t$ is even and $n_i\ge 2$ for every $i\in [t];$
\item $t$ is odd, $n_i\ge 2$ for every $i\in [t]$ and at least two distinct subdigraphs $H_i$ have arcs;
\item  $t$ is odd and $n_i\ge 3$ for every $i\in [t]$ apart from one $i$ for which $n_i\ge 2$.
\end{itemize}
\item[(c)] $T$ and all $H_i$ are strong digraphs with at least two vertices.
\end{description}
\end{thm}

They used Theorem~\ref{sss2} to prove the following characterization for semicomplete compositions $T[H_1,\dots ,H_t]$ when each $H_i$ has at least two vertices.

\begin{thm}\label{sss3}\cite{Sun-Gutin-Ai}
Let $Q=T[H_1,\dots ,H_t]$ be a strong semicomplete composition such that
each $H_i$ has at least two vertices for $i\in [t]$.
Then $Q$ has a strong arc decomposition if and only if $Q$ is not isomorphic to one of the following three digraphs:
$\overrightarrow{C}_3[\overline{K}_2,\overline{K}_2,\overline{K}_2]$, $\overrightarrow{C}_3[\overrightarrow{P_2},\overline{K}_2,\overline{K}_2]$,
$\overrightarrow{C}_3[\overline{K}_2,\overline{K}_2,\overline{K}_3].$
\end{thm}

Bang-Jensen, Gutin and Yeo solved an open problem in \cite{Sun-Gutin-Ai} by obtaining a characterization of all semicomplete compositions with a strong arc decomposition. Note that a digraph with a strong arc decomposition is 2-arc-strong.

\begin{thm}\label{sss4}\cite{Bang-Jensen-Gutin-Yeo}
Let $Q=T[H_1,\dots ,H_t]$ be a semicomplete composition.
Then $Q$ has a strong arc decomposition if and only if $Q$ is 2-arc-strong  and is not isomorphic to one of the following four digraphs: $S_4$,
$\overrightarrow{C}_3[\overline{K}_2,\overline{K}_2,\overline{K}_2]$, $\overrightarrow{C}_3[\overrightarrow{P_2},\overline{K}_2,\overline{K}_2]$,
$\overrightarrow{C}_3[\overline{K}_2,\overline{K}_2,\overline{K}_3].$
\end{thm}

Among their argument, the main technical result is the following theorem.

\begin{thm}\label{sss4-1}\cite{Bang-Jensen-Gutin-Yeo}
Let $Q=T[\overline{K}_{n_1}, \dots,\overline{K}_{n_t}]$ be an extended semicomplete digraph where $n_i\leq 2$ for $i\in [t]$. If $Q$ is 2-arc-strong, then $Q$ has a strong arc decomposition if and only if $Q$ is not isomorphic to one of the following four digraphs: $S_4$,
$\overrightarrow{C}_3[\overline{K}_2,\overline{K}_2,\overline{K}_2]$, $\overrightarrow{C}_3[\overrightarrow{P_2},\overline{K}_2,\overline{K}_2]$.
\end{thm}

Theorems~\ref{intro01} and~\ref{sss4} imply a characterization of quasi-transitive digraphs with a strong arc decomposition (this solves another open question in \cite{Sun-Gutin-Ai}).

\begin{thm}\label{sss5}\cite{Bang-Jensen-Gutin-Yeo}
A quasi-transitive digraph $D$ has a strong arc decomposition if and only if $D$ is 2-arc-strong and is not isomorphic to one of the following four digraphs: $S_4$,
$\overrightarrow{C}_3[\overline{K}_2,\overline{K}_2,\overline{K}_2]$, $\overrightarrow{C}_3[\overrightarrow{P_2},\overline{K}_2,\overline{K}_2]$,
$\overrightarrow{C}_3[\overline{K}_2,\overline{K}_2,\overline{K}_3].$
\end{thm}

All proofs in \cite{Bang-Jensen-Gutin-Yeo} are constructive and can be turned into polynomial algorithms for finding strong arc decompositions. Thus, the problem of finding a strong arc decomposition in a semicomplete composition, which has one, admits a polynomial time algorithm.

Recall that strong semicomplete compositions generalize both strong
semicomplete digraphs and strong quasi-transitive digraphs. However, they
do not generalize locally semicomplete digraphs and their generalizations in-
and out-locally semicomplete digraphs. While there is a characterization of locally semicomplete digraphs having a strong arc decomposition \cite{Bang-Jensen-Huang}, no such a characterization is
known for locally in-semicomplete digraphs and it would be interesting to
obtain such a characterization or at least establish the complexity of deciding whether an locally in-semicomplete digraph has a strong arc decomposition.

\begin{op}\label{op7}\cite{Bang-Jensen-Gutin-Yeo}
Can we decide in polynomial time whether a given locally in-semicomplete digraph has a strong arc decomposition?
\end{op}

\begin{op}\label{op8}\cite{Bang-Jensen-Gutin-Yeo}
Characterize locally in-semicomplete digraphs with a strong arc decomposition.
\end{op}

The following two types of strong subgraph packing problems could be seen as natural extensions of the problem of strong arc decomposition. Let $D=(V(D),A(D))$ be a digraph of order $n$, $S\subseteq V$ a $k$-subset of $V(D)$ and $2\le k\leq n$. A strong subgraph $H$ of $D$ is called an {\em $S$-strong
subgraph} if $S\subseteq V(H)$. Two $S$-strong subgraphs are said to be {\em arc-disjoint} if they have no common arc. Furthermore, two
arc-disjoint $S$-strong subgraphs are said {\em internally disjoint}
if the set of common vertices of them is exactly $S$.
The input of {\sc Arc-disjoint strong subgraph packing} (ASSP) consists of a digraph $D$ and a subset of vertices $S\subseteq V(D)$, the goal is to find a largest collection of arc-disjoint
$S$-strong subgraphs. Similarly, the input of {\sc Internally-disjoint strong subgraph packing} (ISSP) consists of a digraph $D$ and a subset of vertices $S\subseteq V(D)$, and the goal is to find a largest collection of internally
disjoint $S$-strong subgraphs.

Some results on the above strong subgraph packing problems and related topics have been obtained in \cite{Sun-Gutin2, Sun-Gutin-Yeo-Zhang, Sun-Jin, Sun-Zhang} (also can be found in a recent survey \cite{Sun-Gutin3}). Especially, Sun, Gutin and Zhang gave two sufficient conditions for a digraph composition to have at least $n_0$ arc-disjoint $S$-strong subgraphs for any $S\subseteq V(Q)$ with $2\leq |S|\leq |V(Q)|$. Recall that $n_0=\min\{n_i\mid i\in [t]\}$.

\begin{thm}\label{sss6}\cite{Sun-Zhang}
The digraph composition $Q=T[H_1,\dots, H_t]$ has at least $n_0$ arc-disjoint $S$-strong subgraphs for any $S\subseteq V(Q)$ with $2\leq |S|\leq |V(Q)|$, if one of the following conditions holds:
\begin{description}
\item[(a) ] $D$ is a strong symmetric digraph;
\item[(b) ] $D$ is a strong semicomplete digraph and $Q\not\in
\mathcal{Q}_0$, where \\ $\mathcal{Q}_0=\{\overrightarrow{C}_3[\overline{K}_2,\overline{K}_2,\overline{K}_2], \overrightarrow{C}_3[\overrightarrow{P}_2,\overline{K}_2,\overline{K}_2], \overrightarrow{C}_3[\overline{K}_2,\overline{K}_2,\overline{K}_3]\}$.
\end{description}
Moreover, these strong subgraphs can be found within time
complexity $O(n^4)$, where $n$ is the order of $Q$.
\end{thm}

By Theorems~\ref{intro01} and~\ref{sss6}, we directly have:

\begin{cor}\label{sss7}\cite{Sun-Zhang}
Let $Q\not\in \mathcal{Q}_0$ be a strong quasi-transitive digraph. We can in polynomial time find at least $n_0$ arc-disjoint $S$-strong subgraphs in $Q$ for
any $S\subseteq V(Q)$ with $2\leq |S|\leq |V(Q)|$.
\end{cor}

\vskip 1cm

\noindent {\bf Acknowledgement.} Yuefang Sun was supported by Zhejiang Provincial Natural Science Foundation under Grant No.  LY23A010011 and Yongjiang Talent Introduction Programme of Ningbo under Grant No. 2021B-011-G.

\end{document}